\documentclass[10pt]{amsart}
\usepackage{amsmath,amsthm,amsfonts,amssymb,mathrsfs,verbatim,cite,graphicx,hyperref, eucal}

\def\noi{\noindent}

\newtheorem{Thm}{Theorem}[section]
\newtheorem{Def}[Thm]{Definition}
\newtheorem{Lm}[Thm]{Lemma}
\newtheorem{Prop}[Thm]{Proposition}
\newtheorem{Cor}[Thm]{Corollary}

\newtheorem{state}{Theorem}

\theoremstyle{definition}
\newtheorem*{ack}{Acknowledgements}

\numberwithin{equation}{section}

\newlength{\cellsize}
\newcommand\tableau[1]{
\vcenter{
\let\\=\cr
\baselineskip=-16000pt
\lineskiplimit=16000pt
\lineskip=0pt
\halign{&\tableaucell{##}\cr#1\crcr}}}

\newcommand{\tableaucell}[1]{{%
\def \arg{#1}\def \void{}%
\ifx \void \arg
\vbox to \cellsize{\vfil \hrule width \cellsize height 0pt}%
\else
\unitlength=\cellsize
\begin{picture}(1,1)
\put(0,0){\makebox(1,1){$#1$}}
\put(0,0){\line(1,0){1}}
\put(0,1){\line(1,0){1}}
\put(0,0){\line(0,1){1}}
\put(1,0){\line(0,1){1}}
\end{picture}%
\fi}}

\def\cal{\mathcal}
\def\Bbb{\mathbb}
\def\mf{\mathfrak}

\def\<{\langle}
\def\>{\rangle}

\def\a{\alpha}
\def\b{\beta}

\def\th{\theta}
\def\l{\lambda}
\def\L{\Lambda}

\def\e{{\varepsilon}}

\def\Re{\Bbb R}
\def\F{\Bbb F}
\def\C{\Bbb C}
\def\Z{\Bbb Z}

\def\Q{\Bbb Q}

\def\Scal{\cal S}
\def\T{\cal T}
\def\P{\cal P}
\def\H{\cal H}
\def\Qu{\mf Q}

\def\h{\mf h}
\def\g{\mf g}

\def\ag{{\bf a}}
\def\nght{{\bf ht}}
\def\mm{{\cal M}}
\def\m{{\mf m}}
\def\QQ{\cal Q}

\DeclareMathOperator{\hght}{ht}
\DeclareMathOperator{\wt}{wt}
\DeclareMathOperator{\Nght}{Ht}
\DeclareMathOperator{\qwt}{qwt}
\DeclareMathOperator{\Des}{Des}

\begin{document}
\title[]
{Generalized exponents of small representations. I.}
\author{Bogdan Ion}
\thanks{E-mail address: {\tt bion@pitt.edu}}
\date{April 6, 2009}
\address{Department of Mathematics, University of Pittsburgh, Pittsburgh, PA 15260}
\address{Algebra and Number Theory Research Center, Faculty of Mathematics and Computer Science, University of Bucharest, 14 Academiei St., Bucharest, Romania}

\maketitle


\section*{Introduction} 

Kostka polynomials, also known as $t$--weight multiplicities, are ubiquitous in representation theory, geometry, and combinatorics but their structure is far from being unraveled. For the most part, their importance can be attributed to the fact that they are maximal parabolic Kazhdan-Lusztig polynomials for affine Weyl groups. 

Let $\g$ be a complex simple Lie algebra of rank $n$ and denote by $G$ its adjoint group. The Kostka polynomials $m_{\l\mu}(t)$ for $G$ are polynomials in $t$ indexed by pairs of dominant weights of $G$. Denote by $V_\l$ the irreducible representation of $G$ with highest weight $\l$.  The integers $m_{\l\mu}(1)$ are the $\mu$-weight multiplicities of $V_\l$. In general, there are two ways to compute the Kostka polynomials: first, there is an algorithm due to Kazhdan and Lusztig, and second, a formula due independently to Hesselink and Peterson which  computes them in terms of a $t$--analogue of the Kostant partition function. Both approaches become impractical quite quickly, even for computers.

It is possible to obtain some qualitative information on $m_{\l\mu}(t)$. For example, being Kazhdan-Lusztig polynomials, they have non-negative integer coefficients, a fact that can also be seen from the interpretation of their coefficients as dimension jumps of the principal filtration of the $\mu$--weight space of the irreducible representation with highest weight $\l$ (see  \cite{brylinski}). From either description it is  rather difficult to extract any explicit information. 

Since Kostka polynomials have non-negative integer coefficients, many specialists favor the point of view that the best formulas would have to manifest this fact but non-negative formulas which are valid across all Dynkin types are entirely missing from the literature save for the Shapiro-Steinberg procedure for computing the classical exponents of $G$. 
The only other non-negative formula holds only in type $A$. This remarkable formula, due to Lascoux and Sch\" utzenberger \cite{lascoux} (see also \cite{butler}),  expresses each Kostka polynomial as a generating function of the charge statistic over a certain set of tableaux.  Most likely, it is the intricate combinatorics of the charge statistic that foiled all efforts to reformulate the result and its proof in a Lie theoretic manner and extend this formula to other types. 

For one particular case of Kostka polynomials the problem is even older. The polynomials $m_{\l,0}(t)$ (which in this paper will be denoted by $E(V_\l)$) are the  generalized exponents defined by Kostant \cite{kostant}.  
As explained in \cite[$\S$ 2.4]{butler}  crucial insight in the structural properties of  the charge statistic is offered by the charge of  {\sl standard Young tableaux}.
In Lie theoretic terms, this corresponds to understanding the generalized exponents of the so-called first layer representations of $SL_{n+1}(\C)$. Small representations of $SL_{n+1}(\C)$ are either first layer representations or their contragredients and contragredient representations have the same generalized exponents. Therefore, we can consider that the above mentioned formula as corresponding to generalized exponents of small representations. It is also important to mention that in type $A$, thanks to a stability property of Schur functions and Hall-Littlewood functions, any $t$-weight multiplicity equals a $t$-weight multiplicity indexed by small dominant weights but possibly for a larger rank root system.

This is the first paper in a sequence devoted to giving a manifestly non-negative formula for the generalized exponents of small representations {\sl in all types}. The restriction to small representations is not arbitrary but, as we hope to explain in future publications, a crucial stepping stone toward the general situation. The main part of this paper is a complete treatment of the type $A$ case. 

There are several reasons for considering this case separately. First, the argument here should serve as a blueprint for the general argument. Second, this case received more attention than the general one and the ideas and concepts involved in the proof require some discussion vis-\` a-vis what was previously known. Third, although the proof is uniform across types and the argument presented here is rather a specialization of the general argument, certain features specific only to type $A$ were making the general argument unnecessarily complicated and, for clarity reasons, it is easier to treat this case separately. 

The fact that makes everything possible is the computation of the Fourier coefficients of the Cherednik kernel (Theorem \ref{fourierthm}). For this and other finite root systems these coefficients encode the combinatorics of minimal expressions of a small weight (see \cite{ion-geII}). Our main result (Theorem \ref{thmquasi})  is the following. See \eqref{csp} for the definition of the (degenerate) Cherednik scalar product $ \<\cdot,\cdot\>_t$, Definition \ref{quasiorbit} for the concept of  quasi-dominant weight $\l$,  and Definition \ref{quasifundef} for the concept of fundamental quasisymmetric function $\QQ_\l$. The first layer representations in type $A_n$ are those parametrized by partitions of $n+1$. 
\begin{state}\label{state1}
 Let $\l$ be a first layer quasi-dominant weight. Then 
$$
 \<1,\QQ_\l\>_t=t^{\hght(\l)}
$$
\end{state}
The (non-negative integer) coefficients of the expansion of the character of an irreducible small representation $V_\l$ in terms of quasisymmetric functions are called quasi-weight multiplicities $q_{\l\mu}$ and the set of quasi-dominant weights $\mu$ for which $q_{\l\mu}$ is non-zero is denoted by $\qwt(\l)$. As an immediate consequence of Theorem \ref{state1} we obtain 
\begin{state} \label{state2}
Let $\l$ be a first layer dominant weight. Then,
$$
E(V_\l)=\sum_{\mu\in{\rm qwt}(\l)} q_{\l\mu}t^{\hght(\mu)}
$$
\end{state}

As far as I know, the set ${\rm qwt}(\l)$ and the quasi-weight multiplicities were not studied for other root systems. Sufficient information about them can be obtained so that the above formula becomes entirely explicit. In our case this information is already available thanks to work of Gessel. Gessel's result is simply  a consequence of a very basic result of Stanley on $P$-partitions. Next we will describe the outcome. 

As mentioned above, a first layer dominant weight $\l$ in type $A_n$ can be thought of as a partition of $n+1$. Let us denote by $[n]$ the set of integers from $1$ to $n$.
The set of standard Young tableaux of shape $\l$ is denoted by $SYT(\l)$. 
For a fixed standard Young tableau $\T$, the set $\Des(\T)$ of descents consists of the elements $i$ of $[n]$ such that $i+1$ appears in the tableau $\T$ in a row strictly below  $i$. The set of non-descents ${}^c\!\Des(\T)$ is the complement of $\Des(\T)$ inside $[n]$. For every standard Young tableau $\T$ of shape $\l$ one can associate a quasi-weight of $V_\l$ and the height of $\T$ is defined to be the height of this quasi-weight. The definition can be made entirely explicit
$$
\hght(\T)=\sum_{i \in {}^c\!\Des(\T)} (n+1-i)
$$
The integer $\hght(\T)$ turns out to be nothing else but the charge of $\T$. Theorem \ref{state2} now reads 
\begin{state}\label{state3}
Let $\l$ be a first layer dominant weight. Then,
$$
E(V_\l)=\sum_{\T\in{\rm SYT}(\l)} t^{\hght(\T)}
$$
\end{state}
Therefore, Theorem \ref{state3} turns out to be nothing else but the  Lascoux--Sch\" utzenberger charge formula for standard Young tableaux. Our approach to the  charge formula is neither the shortest nor the most simple.  In fact a proof of a statement dual to Theorem \ref{state3} is contained in two examples in \cite{macbook}: Chapter III, $\S$ 6, Example 2 and Chapter I, $\S$ 5, Example 14. It has however some virtues: it offers a Lie theoretic formulation of  charge  as  height of certain weights of $V_\l$ and, more importantly, it offers a proof of the charge formula whose ingredients are available for any root system. The analogues of the results in Section \ref{I-fourier} are contained in \cite{ion-geII} and the analogues of the results in Section \ref{quasifunctions} are the subject of \cite{ion-geIII}. For another approach  that uses concepts available for all root systems see \cite{stembridge, stembridge2}. 

The paper also contains  some intermediate results such as Theorem \ref{firstformulacor} and Theorem \ref{secondformulathm} which compute generalized exponents in terms of weight multiplicities and some combinatorial data such as heights and aggregate vectors of the weights of $V_\l$. These formulas are combinatorially explicit but not manifestly non-negative. 

 \begin{ack}  
 Work supported in part by the NSF grant DMS--0536962 and  the  CNCSIS grant nr. 24/28.09.07 (Groups, quantum groups, corings, and representation theory).  Part of this work was carried out  while visiting the Hausdorff Research Institute for Mathematics, Bonn (June-July 2007) and the Max Planck Institute for Mathematics (June-July 2008);  I thank both institutes for their hospitality and support. I thank J. R. Stembridge for drawing my attention to \cite[I, $\S$ 5, Ex. 14]{macbook} and \cite[III, $\S$ 6, Ex. 2]{macbook}.
 \end{ack}

\section{Preliminaries}


\subsection{} We denote the integers by $\Z$. The word positive, respectively negative,  integer will refer to strictly positive, respectively negative, integers. We use the term non-negative, respectively non-positive, to refer to the set of positive integers and zero, respectively to the set of negative integers and zero.

If $N$ is a positive integer we denote by $[N]$ the set $\{1,2,\dots,N\}$. Also,  $[0]$ refers to the empty set.
If $S$ a set then $\P(S)$ will denote the power set of $S$.

\subsection{}

Let $\g$ be a complex simple Lie algebra of rank $n$ and denote by
$G$ its adjoint group. Let $\h$ and $\mf b$ be  a Cartan
subalgebra respectively a Borel subalgebra of $\g$ such that
$\h\subset {\mf b}$, fixed once and for all. The maximal torus of
$G$ corresponding to $\h$ is denoted by $H$. We have $H=TA$ where
$T$ is a compact torus and $A$ is a real split torus. The volume
one Haar measure on $T$ is denoted by $ds$.

Let $R\subset \h^*$ be the set of roots of $\g$ with respect to
$\h$, let $R^+$ be the set of roots of $\mf b$ with respect to
$\h$ and denote $R^-:=-R^+$. Of course, $R=R^+\cup R^-$; the roots
in $R^+$ are called positive and those in $R^-$ negative. The set
of positive simple roots is denoted by $\{\a_1,\dots,\a_n \}$. We
know that the roots in $R$ have at most two distinct lengths. We
will use the notation $R_s$ and $R_\ell$ to refer respectively to
the short roots and the long roots in $R$. The dominant
element of $R_s$ is denoted by $\th_s$ and the dominant element of
$R_\ell$ is denoted by $\th_\ell$.

Any element $\a$ of $R$ can be written uniquely as a sum of simple
roots $\sum_{i=1}^n a_i\a_i$. The height of the root $\a$ is
defined to be $$\hght(\a)=\sum_{i=1}^n a_i$$ The root of $R$ with
has the largest height is denoted by $\th$. By the above
convention, if $R$ is simply laced then $\th=\th_s$ and if $R$ is
not simply laced then $\th=\th_\ell$.

Denote by $r$ the maximal number of laces in the Dynkin diagram
associated to $\g$. There is a canonical positive definite
bilinear form $(\cdot,\cdot)$ on $\h^*_\Re$ (the real vector space
spanned by the roots) normalized such that $(\a,\a)=2$ for long
roots and $(\a,\a)=2/r$ for short roots.  For any root $\a$ define
$\a^\vee:=2\a/(\a,\a)$. We know from the axioms of a root system
that $(\a,\b^\vee)$ is an integer for any roots $\a$ and $\b$. In
fact, the only possible values for $|(\a,\b^\vee)|$ are $0$, $1$
or $2$ if the length of $\a$ does not exceed the length of $\b$
(the value $2$ is attained only if $\a=\pm\b$) and $0$, $r$ if the
length of $\a$ is strictly larger than the length of $\b$. It is a well-known fact that if $R$ is not simply laced then 
$$
(\th_\ell,\th_s)=1
$$

Define $\rho=\frac{1}{2}\sum_{\a\in R}\a^\vee$. With this notation
the height of any root $\a$ can be written as
$\hght(\a)=(\a,\rho)$. The weight lattice of $G$ is $Q$,  the integral span
of the simple roots; we will use the word weight to refer to an element of $Q$. For a weight $\l$  define its
height as $\hght(\l)=(\l,\rho)$.

For any root $\a$ consider the reflection of the Euclidean space
$\h^*_\Re$ given by $$s_\a(x)=x-(x,\a^\vee)\a .$$ The Weyl group
$W$ of the root system $R$ is the subgroup of ${\rm GL}(\h^*_\Re)$
generated by the reflections $s_\a$, for all roots $\a$ (the
simple reflections $s_i:=s_{\a_i}$, $1\leq i\leq n$, are enough).
The scalar product on $\h^*_\Re$ is equivariant with respect to
the action of $W$.

For any $w$ in $W$ denote by $\ell(w)$  the length of a reduced (i.e.
shortest) decomposition of $w$ in terms of simple reflections. The element $w_\circ$ is the unique maximal
length element in $W$.


\subsection{}

For an element $\l$ of the root lattice  we denote by $e^\l$ the
corresponding character of the compact torus $T$. The trivial
character $e^0$ will be also denoted by $1$. Let $\Z[Q]$ be the
$\Z$--algebra spanned by all such elements (the group algebra of
the lattice $Q$). Note that the multiplication is given by
$e^\l\cdot e^\mu=e^{\l+\mu}$.

The subalgebra of $\Z[Q]$ consisting of $W$--invariant elements is
denoted by $\Z[Q]^W$.  The irreducible finite dimensional
representations of $G$ are parameterized by the dominant weights. For a dominant $\l$ we denote by $\chi_\l$
the character of the corresponding irreducible representation of
$G$. Restricting the characters to $T$ we will regard them as
elements of $\Z[Q]$. A basis of $\Z[Q]^W$ is then given by the all
the irreducible characters $\chi_\l$ of $G$.


\subsection{}

Assume that $t$ is a complex number of small absolute value and
let
$$C(t)=\prod_{\a\in
{R^+}}\frac{1-e^{\a}}{1-te^{\a}}$$ Since  $t$ is small the infinite product is absolutely convergent
and $C(t)$ should be seen as a continuous function on the
torus $T$. The function $C(t)$ is the specialization at $q=0$ of the so--called Cherednik kernel $
C(q,t)$ (see, for example \cite[(3.9)]{ion-ge}). Of course,  $t$ can be considered as a formal variable, in which case we have to work over the field
$\F:=\Q(t)$.

Consider the $\F$-linear involution of $\F[Q]$ given by $\overline{e^\l}=e^{-\l}$. The following pairing  is a non--degenerate scalar product on
$\F[Q]$
\begin{equation}\label{csp}
\<f,g\>_{t}:=\int_T f\overline{g} C(t)ds
\end{equation}
This is simply the specialization at $q=0$ of the usual Cherednik scalar product (see, for example, \cite[(3.11)]{ion-ge}).

Let $\vartheta$ be an automorphism of the Dynkin diagram of $\g$. The action of $\vartheta$ on the simple roots can be extended linearly to an action on $Q$ and then to an $\F$-linear action on $\F[Q]$. Since this action leaves invariant the set of positive roots, $C(t)$ is fixed by $\vartheta$. As a consequence, $\vartheta$ is unitary for the above scalar product
\begin{equation}\label{varsigmaunitary}
\<\vartheta(f),\vartheta(g)\>_{t}=\<f,g\>_{t}
\end{equation}
For the purposes on this paper it is enough to consider the automorphism  $\vartheta_\circ$ of the Dynkin diagram that sends $e^\l$ to
$e^{-w_\circ(\l)}$. 


\subsection{} 

Let $\chi_\l$ denote the character of $V_\l$, the irreducible representation
of $G$ with highest weight $\l$ and let ${\rm wt}(\l)$ be the set
of weights of $V_\l$.  As it is well-known, ${\rm wt}(\l)$ consists of the elements of the root lattice which are contained in the convex hull of the Weyl group orbit of $\l$.
For any $\gamma\in {\rm wt}(\l)$ we use
 $m_{\l \gamma}$ to refer the weight multiplicity of $\gamma$ in
$V_\l$.

The zero weight space of $V_\l$, denoted by $V_\l(0)$, will also be important for us. The torus $T$ acts trivially on $V_\l(0)$ and therefore its normalizer $N_G(T)$ in $G$ acts on the zero weight space. Keeping in mind that the Weyl group $W$ is isomorphic to $N_G(T)/T$ we obtain that $V_\l(0)$ is a $W$ representation.


\subsection{}\label{kostantfactorization}
Consider now $t$ as a formal variable. The graded torus
character of $\S(\g)$, the algebra of complex valued polynomial
functions on $\g$, is easily seen to be
$$\frac{1}{(1-t)^n}\prod_{\a\in R}\frac{1}{1-te^\a}$$Kostant \cite{kostant} studied the action of $G$ on the symmetric algebra ${\Scal}(\g)$ of complex valued polynomial functions on $\g$. One of his fundamental results is the following.  Denote by $\Scal(\g)^G$  the subring of $G$--invariant polynomials on $\g$. 
Then, $\Scal(\g)$ is free as an $\Scal(\g)^G$--module and it is generated by $\H(\g)$, the space of $G$--harmonic polynomials on $\g$ (the polynomials annihilated by all $G$--invariant differential operators with constant complex coefficients and no constant term). In other words,
$$\Scal(\g)=\Scal(\g)^G\otimes \H(\g)$$ The space of harmonic polynomials as a graded, locally finite representation of $G$ is itself of considerable interest, being part of rich theory at the intersection of geometry and representation-theory. The space $\H(\g)$, which in fact inherits an algebra structure from the symmetric algebra, can be regarded as the
 the ring of regular functions on the cone of nilpotent elements in $\g$. The $G$--module structure of rings of regular functions on closures of nilpotent orbits lie at the heart of many questions about primitive ideals of enveloping algebras, associated varieties and characteristic cycles (see, for example, \cite{vogan}) and explicit results are desirable not only in the context of complex semisimple Lie algebras but much more generally. From this point of view the simplest situation is the one we will be concerned with here: understanding the $G$-module structure of $\H(\g)$.

Denote by $\H^i(\g)$ the degree $i$ component of $\H(\g)$, and denote by $V_\l$ the irreducible
representation of $G$ with highest weight $\l$.  The
graded multiplicity of $V_\l$ in $\H(\g)$
$$
E(V_\l):=\sum_{i\geq 0} {\rm dim_\C}\left({\rm
Hom}_G(V_\l,\H^i(\g))\right)t^i
$$
is a polynomial with non-negative integer coefficients. The multiplicity of $V_\l$ inside $\H(\g)$ is obtained from substituting $1$ for $t$ in  the above formula. This multiplicity equals $v_\l$ the dimension of the zero weight space of $V_\l$.

It is therefore possible to write  $$E(V_\l)=\sum_{i=1}^{v_\l} t^{e_i(\l)}$$
such that $e_1(\l)\leq e_2(\l)\leq \cdots \leq e_{v_\l}(\l)$.  The
positive integers $e_i(\l)$ are important invariants of the representation $V_\l$, first defined and studied by Kostant in \cite{kostant} who called them the generalized
exponents of $V_\l$. As Kostant explained in his work, the terminology is justified by the fact that the generalized exponents of the adjoint representation of $G$ are in fact the classical exponents $e_1\leq \cdots \leq e_n$  of $G$, which are extremely basic invariants that appear in many important contexts such as  the topology of Lie groups and the invariant theory of real reflection groups.

The  graded multiplicity of
$V_\l$ inside $\H(\g)$ is denoted by $E(V_\l)$. Let us recall the following observation from \cite[(3.18), (3.13)]{ion-ge}
\begin{equation}
E(V_\l)=\<1,\chi_\l \>_{t}
\end{equation}
As it is well-known, and also immediately follows from  \eqref{varsigmaunitary},  
\begin{equation}\label{contragredient}
E(V_{\vartheta(\l)})=E(V_\l)
\end{equation}
for any $\vartheta$ an automorphism of the Dynkin diagram of $\g$. In particular, the above equality for the involution $\vartheta_\circ$ implies that 
the relevant graded multiplicities for a representation and its contragredient are equal.

\subsection{}

The elements $$c_\mu(t)=\<1,e^\mu\>_{t}$$   are known to be elements of $\F$ and, as pointed out in \cite{ion-ge}, have a major role to play in the computation of  $E(V_\l)$. 
For any continuous function $f$ on the torus $T$, its Fourier
coefficients are parametrized by weights and are given by
$$
f_\l:=\int_T fe^{-\l}ds .
$$
From this point of view $c_\l(t)$ are the Fourier coefficients of  $C(t)$ (regarded as a continuous function on $T$).

In \cite[section 5.1]{ion-ge} was described a linear system which has as unique solution the Fourier coefficients  of the Cherednik kernel. 
To streamline subsequent computations we consider here an equivalent system. Let us start by recording the specialization at $q=0$ of  \cite[(5.4) and (5.8)]{ion-ge}
\begin{subequations}
\begin{align} \label{eq1}
c_{s_i(\l)}(t)-t^{-1}c_\l(t)&=(t^{-1}-1)\left( c_{\l-\a_i}(t)+\cdots
+c_{\l-(k-1)\a_i}(t) \right)\\
\intertext{for any weight $\l$ and $\a_i$ a simple root such that $(\l,\a_i^\vee)=k>0$.
 Also, }
\label{eq2}
c_{s_\th(\l)}(t)&=0
\end{align}
\end{subequations}
if $\l$ is dominant. 
\begin{Lm} Let $\l$ be a weight and let $\a_i$ be a simple root such that $(\l,\a_i^\vee)>0$.
Then,  \begin{equation}\label{eq3}
c_{s_i(\l)}(t)-t^{-1}c_\l(t)=- c_{\l-\a_i}(t)+
t^{-1}c_{s_i(\l)+\a_i}(t) 
\end{equation}
\end{Lm}
\begin{proof} If $(\l,\a_i^\vee)\leq 2$ then it is easy to check that  equation \eqref{eq3} is exactly 
 \eqref{eq1}. If  $(\l,\a_i^\vee)> 2$ then $(\l-\a_i,\a_i^\vee)=(\l,\a_i^\vee)-2>0$ and \eqref{eq3} is obtained by subtracting the equations \eqref{eq1} for $\l$ and $\l-\a_i$.
\end{proof}

In fact, the equations \eqref{eq2} and \eqref{eq3}  determine the coefficients $c_\l(t)$.  The reason for that  is essentially the one from \cite[Theorem 5.1]{ion-ge} for 
$q=0$, but in order to stress and to clarify one important aspect  that will be used in this paper we recall here the argument.

Fix a  dominant weight $\l_+$ and consider the following
homogeneous linear system:  the unknowns
are $x_\l$  for all $\l\in W\l_+$ and the equations are
\begin{subequations}\label{eqs4-5}
\begin{alignat}{2}\label{eq4}
& x_{s_i(\l)}-t^{-1}x_\l=0\quad   &\text{if  $(\l,\a_i^\vee)>0$} \\ \label{eq5}
& x_{s_\th(\l_+)} = 0               &
\end{alignat}
\end{subequations}
It is clear that  $x_\l=0$ for all $\l\in W\l_+$  is the unique solution of this system.

\begin{Thm}\label{newsystem}
Let $\l_+$ be an arbitrary dominant weight. Consider the finite linear system ${\rm Sys}(\l_+)$ with unknowns $y_\mu$  indexed by  $\mu\in {\rm wt}(\l_+)$ and equations
\begin{subequations}
\begin{alignat}{2}\label{eq6}
&y_{s_i(\mu)}-t^{-1}y_\mu=- y_{\mu-\a_i}+ t^{-1}y_{s_i(\mu)+\a_i} &\quad & \text{if  $(\mu,\a_i^\vee)>0$} \\ \label{eq7}
&y_{s_\th(\mu)} =  0                                                                                  &           & \text{if $\mu\neq 0$ dominant}\\ \label{eq8}
&y_0=1                                                                                                         &         &
\end{alignat}
\end{subequations}
Then, ${\rm Sys}(\l_+)$ has a unique solution. In particular, the infinite system ${\rm Sys}(\infty)$  described by the same equations but with variables $y_\mu$ indexed by all weights $\mu$ has a unique solution.
\end{Thm}
\begin{proof}
We know that the system has at least one solution since this is guaranteed by equations \eqref{eq2} and \eqref{eq3}. 
We prove that the solution is unique by induction on the distance between $\l_+$ and the origin. 

The solution of the system ${\rm Sys}(0)$ is unique by \eqref{eq8}. Assume now that the distance $d$ from $\l_+$ to the origin is strictly positive and that the system ${\rm Sys}(\mu_+)$ has unique solution for any dominant weight that is strictly closer to the origin.  All the elements of ${\rm wt}(\l_+)$ (which are exactly the weights in the convex hull of $W\l_+$) except those in $W\l_+$ have distance to the origin strictly smaller than $d$ and therefore they are uniquely determined by the induction hypothesis. Hence we only have to argue that $y_\l$  for all $y\in W\l_+$ are uniquely determined.  

For this it is enough to show that the system with unknowns  $\mu\in W\l_+$ and equations \eqref{eq6} and \eqref{eq7}  has at most one solution. Note that in this case the right hand side of \eqref{eq6} is already known by the induction hypothesis. The associated homogeneous system  is the one  given by equations \eqref{eqs4-5} which has has a unique solution and this implies our claim. The statement about ${\rm Sys}(\infty)$ is an immediate consequence.
\end{proof}
\begin{Cor}\label{newsystemverification} 
Let $\Gamma\subseteq \h_\Re^*$ be a union of convex hulls of $W$-orbits. The system ${\rm Sys}(\Gamma)$ obtained from ${\rm Sys}(\infty)$ by restricting to variables and equations involving only $\mu$ in ${Q\cap \Gamma}$ has a unique solution.
\end{Cor}
\begin{proof}
The result  follows from the above Theorem keeping in mind that $Q\cap \Gamma$ can be written as a union of sets of the form ${\rm wt}(\l)$ for $\l$ dominant weights.
\end{proof}
This simple result is one of our the main tools: if $\Gamma$ is as above, to prove any conjectural formula for the Fourier coefficients $c_\l(t)$ for $\l$ in $Q\cap\Gamma$   it is enough to show that the formula satisfies the system ${\rm Sys}(\Gamma)$.
\section{Small representations}

\subsection{} Let $\l$ be a dominant weight. The highest weight representation $V_\l$ is called small if and only if 
\begin{equation}\label{smallcondition}
2\theta_s\not \in \wt(\l)
\end{equation}

Any weight that is a weight for some  small representation will be called a small weight. Of course, a weight $\l$ is small if and only if $\l_+$, the unique dominant element in $W\l$, is small, which holds true exactly when $V_{\l_+}$ is a small representation. We will denote by $Q^{\rm sm}$ the set of small weights.

A  parametrization of small dominant weights by their canonical block decomposition is provided in \cite[Theorem 1]{ion-geII}. In type $A$ the small representations are either first layer representations or their contragredients. We discuss this in Section \ref{layer-k}.


\subsection{}\label{broer}  There is a beautiful connection between the representation theory of $G$ and that of $W$ which stresses the special role played by the small representations.  

Kostant's    tensor product decomposition for the symmetric algebra $\Scal(\g)$ mentioned in Section \ref{kostantfactorization} was preceded  by an analogous result of Chevalley regarding the symmetric algebra $\Scal(\h)$. If $\Scal(\h)^W$ denotes the subring of
$W$--invariant polynomials on $\h$
 and
$\H(\h)$ the space of $W$--harmonic polynomials on $\h$
then Chevalley's  result states that $$\Scal(\h)=\Scal(\h)^W\otimes \H(\h)$$ as graded $W$--modules. More information is available in this case: $\H(\h)$ is usually referred to  in the literature as the covariant ring  of the reflection representation of $W$; as an algebra it is isomorphic to the cohomology of the flag variety of $G$ and as a $W$ representation it is nothing else but the regular representation of $W$. Chevalley also related the two invariant algebras: the restriction map sending polynomials on $\g$ to their restriction on $\h$ induces a graded algebra isomorphism
$$
\Scal(\g)^G\to\Scal(\h)^W
$$
In consequence, both algebras have the same Poincar\' e series. 

If $V$ is an irreducible $W$ representation let us denote by $F(V)$ the graded multiplicity of
$V$ inside $\H(\h)$. The polynomial $F(V)$ is known in the literature as the fake degree of $V$; it is, in some sense, closely related to the dimension of a unipotent representation of the group $G$ over a finite field.

Broer \cite{broer} made  the following remarkable observation relating  generalized exponents and  fake degrees.
 \begin{Thm}\label{broerthm} The dominant weight $\l$ is small if and only if 
$E(V_\l)=F(V_\l(0))$.
\end{Thm}
In general, the map that sends a small $G$ representation to the $W$ representation afforded by its zero weight space  neither sends  an irreducible representation to an irreducible representation nor all irreducible $W$ representations arise this way. Both these facts hold true in type $A$. In principle, with the correct understanding of which representations of $W$ arise from zero weight spaces of small representations (some detailed information is available thanks to the work of Reeder \cite{reeder2}), one can use the explicit formulas for fake degrees given by Steinberg \cite{steinberg} (type $A$),  Lusztig \cite[$\S$ 2.4-5]{lusztigfake} and Stembridge \cite[$\S$ 5]{stembridge1} (type $B$, $C$, $D$), and Beynon and Lusztig (type $E$, $F$, $G$) to compute the generalized exponents of small representations. None of these formulas is Lie theoretic in nature and offer little insight on what should be true in general. 


\section{First layer Fourier coefficients in type \texorpdfstring{$A$}{A}} \label{I-fourier}

\subsection{}\label{conventions}  As a general convention, if $S$ is a subset the integers we use the notation
\begin{equation}\label{I-convention}
(1-t^S)
\end{equation}
to refer to $1$ if $S$ is the empty set and to $\prod_{s\in S}(1-t^{\min\{0,s\}})$ otherwise. The product is zero unless $S$ consists of negative integers. We will use the analogue notation for $(t^S-1)$. If ${\bf v}=(v_1,\dots,v_k)$ is a vector in $\Re^k$ and its coordinates in the usual standard basis are integers then we use 
\begin{equation}\label{I-conventionvec}
(1-t^{\bf v})
\end{equation}
to refer to $(1-t^S)$ with $S=\{ v_1,\cdots ,v_k\}$. 

The zero vector in $\Re^k$ will be denoted by $\bf 0$. Note that we suppressed any reference to $k$ from the notation this information being hopefully unambiguous from the context.
If $\bf v$ and $\bf w$ are two vectors in $\Re^k$ we write $$\bf v< w $$ if and only if $v_i<w_i$ for all $1\leq i\leq k$. For any fixed $1\leq i\leq k$ we denote by $\hat{\bf v}^{i}$ the vector in $\Re^{k-1}$ obtained by omitting the $i$-th coordinate from $\bf v$. We use the same notation in the case we need to omit more than one coordinate.


\subsection{}

For the rest of the paper we restrict ourselves to  root systems type $A_n$, $n\geq 1$. Denote by 
 $\{\e_i\}_{1\leq i\leq n+1}$ the standard basis of  $\Re^{n+1}$ and by $(\cdot,\cdot)$ its canonical scalar product. Let $V$ be the subspace of $\Re^{n+1}$ orthogonal to $${\bf 1}:=\e_1+\cdots+\e_{n+1}$$ The root system of the Lie algebra $\mf{sl}(n+1,\C)$ with respect to its Cartan subalgebra $\h$ consisting of diagonal matrices is of type $A_n$. Choose also a Borel subalgebra as the subspace of upper-triangular matrices.
 The Euclidean vector space $\left(\h^*_\Re, (\cdot,\cdot)\right)$ can be identified to $\left(V, (\cdot,\cdot)\right)$. Under this identification  the root lattice $Q$ equals $\Z^{n+1}\cap V$. Moreover,  
 $$R=\{\e_i-\e_j\}_{1\leq i\neq j\leq n+1}\quad \text{and} \quad R^+=\{\e_i-\e_j\}_{1\leq i< j\leq n+1}$$
The simple roots are $$\a_i=\e_i-\e_{i+1}, \quad 1\leq i\leq n$$
and $\th=\e_1-\e_{n+1}$.  The Weyl group is the symmetric group $S_{n+1}$ and its action on $V$ is the usual action permuting the coordinates. 

If we use the notation $$x_i:=e^{\e_i} $$
for all $1\leq i\leq n+1$, the ring $\F[Q]$ becomes  $\F[x_1,\dots, x_{n+1}]/(x_1\cdots x_{n+1}-1)$.
\subsection{}\label{layer-k}  An element $\l=(\l_1,\dots,\l_{n+1})$ of the root lattice is said to be a $k$-th layer weight if its smallest coordinate equals $-k$. For example, the zero weight is the unique 0-th layer weight and the roots are all in the first layer. To simplify some of the later statements we will abuse terminology and consider the zero weight to be also in the first layer. If $\l$ is dominant in the $k$-layer then $\l+k\bf 1$ is a partition of $k(n+1)$.

It is well known that the small dominant weights in type $A_n$ are exactly the first layer dominant weights and their contragredients. For the purpose of computing the generalized exponents of small representations \eqref{contragredient} allows us to restrict to representations whose highest weight is in the first layer. The set of first layer weights  will be denoted by $Q^{(1)}$.

The irreducible representations of $S_{n+1}$ are parametrized by partitions of $n+1$. As already mentioned in Section \ref{broer} all the irreducible representations of $S_{n+1}$ can be realized on zero weight spaces of first layer representations of $\mf{sl}(n+1,\C)$. If $\l$ is  a first layer dominant weight then the zero weight space $V_\l(0)$ affords the irreducible representation of $S_{n+1}$ indexed by the partition dual to $\l+\bf 1$. We refer to \cite{gutkin, kostant2} for  details.

\begin{Def}
Let $\l$ be a first layer weight.  The \emph{length} and \emph{co-length} of $\l$ are defined to be the number of non-negative coordinates of $\l$ and, respectively, 
the number of negative coordinates of $\l$ (which are necessarily equal to $-1$). 
\end{Def}
Note that zero weight in $A_n$ has length $n+1$ and co-length $0$.
It is clear from the definition that the length and co-length are in fact invariants of the Weyl group orbit of $\l$. We will use the notation $\ell(\l)$ and $\ell^*(\l)$ to refer to the length and co-length of $\l$, respectively. If $\l$ is dominant the terminology is consistent to the usual terminology used for partitions. Indeed,  $\ell(\l)$ equals the length of the partition $\l+\bf 1$  of $n+1$ and $\ell^*(\l)$ equals the difference between $n+1$ and $\ell(\l)$.

\subsection{} Assume that $\l=(\l_1,\cdots,\l_{n+1})$ has co-length $N$. Let $$\l_{i_1}=\l_{i_2}=\cdots =\l_{i_N}=-1, \quad i_1<\cdots <i_N$$ be its negative coordinates. Define
\begin{equation}\label{agregate}
\ag_\l(j):=\sum_{k=i_j}^{n+1} \l_k
\end{equation}

There are some basic inequalities  relating these integers. Keeping in mind that $$\sum_{k=1}^{n+1} \l_k=0$$ we obtain that 
\begin{subequations} 
\begin{align}\label{ineq1}
\ag_\l(1)&\leq 0\\ \label{ineq2}
\ag_\l(i)&\geq \ag_\l({i+1})-1, \quad 1\leq i\leq N-1 \\ \label{ineq3}
\ag_\l(N)&\geq -1 
\end{align}
\end{subequations}

\begin{Def}
Let $\l$ be a first layer weight of co-length $N$. The vector \begin{equation}\label{agregatevec}\ag_\l:=(\ag_\l(1),\ag_\l(2),\cdots,\ag_\l(N))\end{equation} will be called the {aggregate vector} of $\l$. Note that $\ag_0$ is  the empty set.
\end{Def}

\subsection{} We derive next two simple facts which will be used in the proof of the main result of this section. 

\begin{Lm}\label{lemma1} Let $\l$ be a first layer weight of positive co-length. Assume that either $s=1$ or $2 \leq s\leq \ell^*(\l)$ and 
$  \ag_\l({s-1})<0$. Then, $\ag_\l(s)\leq0$.
\end{Lm}
\begin{proof}
If $s=1$ then \eqref{ineq1} assures that the conclusion is satisfied. If $s>1$ then from \eqref{ineq2} we obtain that $$0>\ag_\l({s-1}) \geq  \ag_\l({s})-1$$ which implies the desired equality.
\end{proof}

\begin{Lm}\label{lemma2}
Let $\l$ be a first layer weight  and $\a_i$ a simple root such that $(\l,\a_i^\vee)>0$.  
\begin{enumerate}
\item[(a)]  If $\l_{i+1}\neq -1$, then $\ag_{s_i(\l)}=\ag_{\l-\a_i}=\ag_{s_i(\l)+\a_i}=\ag_\l$.
\item[(b)] If $\l_{i+1}= -1$ and this is the $j$-th negative coordinate of $\l$, then $$\hat{\ag}^j_{s_i(\l)}=\hat{\ag}^j_{\l}\quad \text{and}\quad \ag_{s_i(\l)}(j)=\ag_{\l}(j)+\l_{i}$$ Moreover, if $\l_i>0$ then $\l-\a_i$ and $s_i(\l)+\a_i$ have  co-length  $\ell^*(\l)-1$ and $$\ag_{\l-\a_i}=\ag_{s_i(\l)+\a_i}=\hat{\ag}^j_\l$$
\end{enumerate}
\end{Lm}
\begin{proof} (a) The hypothesis forces $\l_i>\l_{i+1}\geq 0$. Therefore, the $i$--th and $i+1$--st coordinates of $\l$, $s_i(\l)$, $\l-\a_i$, and $s_i(\l)+\a_i$  are all
non-negative and their sum is the same in all four cases. The conclusion is immediate.  

\noi (b) Straightforward verification. 
\end{proof}


\subsection{} Our first result is the following.
\begin{Thm}\label{fourierthm} Let $\l$ be a first layer weight. Then, 
\begin{equation}\label{fourier}
c_\l(t)=t^{\hght(\l)}(1-t^{\ag_\l})
\end{equation}
\end{Thm}
\begin{proof} Let $\Gamma$ be the convex hull of $Q^{(1)}$.  It is easy to check that $\Gamma$ is a convex, $S_{n+1}$-invariant set and that $Q\cap \Gamma=Q^{(1)}$. We check that the proposed formula for $c_\l(t)$  satisfies the system ${\rm Sys}(\Gamma)$.

The equation \eqref{eq8} is clearly satisfied.  Assume now that $\l$ is a first layer weight of co-length $N\geq 1$. The  equation \eqref{eq7} is very easy to check. Indeed, assume that $\l$ is dominant. Then,  $\l_{n+1}=-1$ and $\ag_1({s_\th(\l)})= 0$. By our conventions in Section \ref{conventions} this implies $c_{s_\th(\l)}(t)=0$.

Let us argue now that the proposed formula satisfies equation \eqref{eq6}. Let $\a_i$ be a simple root such that $(\l,\a_i^\vee)=k>0$. Denote by $h$ the height of $\l$.

Assume first that $\ag_{\l-\a_i}\not< \bf 0$. But then Lemma \ref{lemma2} implies that the same is true for $\ag_{s_i(\l)}$, $\ag_{s_i(\l)+\a_i}$, and $\ag_\l$. In consequence, all the terms appearing in equation \eqref{eq6} are zero.

Assume now that $\ag_{\l-\a_i}< \bf 0$. From Lemma \ref{lemma2} we know that either $$\ag_{s_i(\l)}=\ag_{\l-\a_i}=\ag_{s_i(\l)+\a_i}=\ag_\l$$  or,
$\l-\a_i$ and $s_i(\l)+\a_i$ have  co-length  $N-1$ and $$\ag_{\l-\a_i}=\ag_{s_i(\l)+\a_i}=\hat{\ag}^j_\l$$ for some  $j$. 

In the former case \eqref{eq6} is trivially satisfied. In the latter case, remark that by applying Lemma \ref{lemma1} for $\l$ we obtain that $\ag_\l(j)\leq 0$ and by applying it for $s_i(\l)$  we obtain that $\ag_{s_i(\l)}(j)\leq 0$.  The equation \eqref{eq3},  up to the common factor  $(1-t^{\ag_{\l-\a_i}})$,  now reads
$$
t^{h-k}(1-t^{\ag_\l(j)+\l_i})-t^{h-1}(1-t^{\ag_\l(j)})=-t^{h-1}+t^{h-k}
$$
Keeping in mind that in this situation $\l_{i+1}=-1$ and that $k=\l_i-\l_{i+1}$ one can easily check that the above equality is satisfied. 

In conclusion, the proposed expressions satisfy the system specified in Corollary \ref{newsystemverification} and therefore they must be  the relevant Fourier coefficients.
\end{proof}

\begin{Thm}\label{firstformulacor}
Let $\l$ be a first layer dominant weight. Then,
\begin{equation}\label{firstformula}
E(V_\l)=\sum_{{\mu\in\wt(\l)}}m_{\l\mu}t^{\hght(\mu)}(1-t^{\ag_\mu})
\end{equation}
\end{Thm}

\subsection{}
As explained and illustrated in \cite{ion-ge}  one can obtain some combinatorial formulas for generalized exponents as an immediate consequence of formulas like \eqref{firstformula}. 
Let  $\l$ be a  first layer dominant weight and let $i$ be a non-negative integer.
Denote by $h_\l^{\rm even}(i)$ and, respectively $h_\l^{\rm odd}(i)$, the number pairs $(\mu,A)$ where $\mu$ is a  weight of $V_\l$ (counted with multiplicity) such that $\ag_\mu<{\bf 0}$  and  $A$ is a subset of $[\ell^*(\mu)]$ of even, respectively odd,  cardinality such that 
$$
\hght(\mu)+\sum_{a\in A} \ag_\mu(a)=i
$$
The result alluded to is the following
\begin{Thm}\label{secondformulathm}
Let $\l$ be a first layer dominant weight. Then,
\begin{equation}\label{secondformula}
E(V_\l)=\sum_{i=1}^{\hght(\l)} (h_\l^{\rm even}(i)-h_\l^{\rm odd}(i)) t^i
\end{equation}
\end{Thm}
\begin{proof}
Straightforward consequence of \eqref{firstformula}.
\end{proof}
Such formulas can  be used to extract completely explicit information for first layer dominant weights of small co-length. For example, for dominant weights of co-length 1 (i.e. the dominant root) \eqref{secondformula} is exactly the Shapiro-Steinberg procedure for computing classical exponents. The dominant weights of  co-length 2 is treated in \cite{ion-ge}. However, as the co-length grows, the combinatorics of  $h^{\rm even}_\l(i)$ and $h_\l^{\rm odd}(\l)$ becomes quite complicated and a more refined formula is therefore desirable. This will be achieved in the next section. 


\section{Quasisymmetric functions}\label{quasifunctions}


\subsection{}\label{quasifunctions1} Quasisymmetric functions were defined by Gessel \cite{gessel} who studied them in connection to enumeration problems for permutations such as counting permutations with given descent set.  We describe this notion in a fashion that is more suitable when dealing with root systems. We postpone the discussion on the relationship between the two constructions until  Section \ref{quasipoly}.

Let $\l$ be a weight and $s_i$ a simple reflection. Define
\begin{equation}\label{localsimpleaction}
s_i\centerdot\l=\begin{cases}
\l& \text{if $\l_i,\l_{i+1}\geq 0$}\\
s_i(\l)& \text{otherwise}
\end{cases}
\end{equation}
It is straightforward to verify that \eqref{localsimpleaction} extends to an action 
\begin{equation}\label{localaction}
S_{n+1}\times Q\to Q, \quad(w,\l)\mapsto w\centerdot\l
\end{equation}
Since the way an element acts is influenced by local conditions we refer to this action as the local action of $S_{n+1}$ on $Q$ and we call local orbit a $S_{n+1}$ orbit with respect to the local action. It is clear from the definition that the local orbit of a weight is a subset of its usual orbit. In the local orbit of $\l$  there is a unique  element of maximal height: the element in the usual orbit with the negative entries in the decreasing order on the leftmost possible positions.
\begin{Def}\label{quasiorbit} 
A first layer weight $\l$ is said to be quasi-dominant if it is the maximal height  element of its local orbit. The set of first layer quasi-dominant weights will be denoted by $Q^{(1), q+}$.
\end{Def}
Let $\l$ be a  quasi-dominant weight.  The function 
$$\mm_\l:=\sum_{\mu\in S_{n+1}\centerdot\l} e^\mu$$
is called a  quasisymmetric monomial. 

The space of functions in ${\rm span}_\F\{e^\l~|~\l\in Q^{(1)}\}$ which are constant on local orbits is called the space of (first layer) quasisymmetric functions and denoted by $\QQ^{(1)}$. The elements in $\QQ^{(1)}$ will be called quasisymmetric functions. It is clear the first layer quasisymmetric monomials form a basis of the space of quasisymmetric functions.

\subsection{}

Assume that $\l$ is a first layer quasi-dominant weight of co-length $N$. There is a unique positive root  (let us call it $\b_N$) such that $\l-\b_N$ has co-length $N-1$ and $\b_N$ is of smallest height with this property. Indeed, if the rightmost positive coordinate of $\l$ is on position $i$ and its leftmost negative coordinate is on position $j$ (in fact $j=n+2-N$) then $\b_N=\e_i-\e_{j}$. 

In fact, more is true: $\l-\b_N$ is itself quasi-dominant. This leads us to the following concept.

\begin{Def}
Let  $\l$ be a first layer quasi-dominant weight of co-length $N$. The canonical expression of $\l$ is defined  inductively to be the expression
$$
\l=\b_1+\cdots+\b_N
$$
where $\b_N$ is the unique positive root  of smallest possible height such that $\l-\b_N$ has co-length $N-1$ and $$\l-\b_N =\b_1+\cdots+\b_{N-1}$$ is the canonical expression of $\l-\b_N$. The canonical expression of the zero weight is by definition  $\l=0$.
\end{Def}
For example the canonical expression of $\l=(0,2,0,1,0,0,-1,-1,-1)$ (which is a first layer quasi-dominant weight for $A_8$ of co-length $3$) is 
$$
\l=(\e_2-\e_9)+(\e_2-\e_8)+(\e_4-\e_7)
$$
It is clear from the definition that heights of the positive roots in the canonical expression of a quasi-dominant weight $\l$ of co-length $N$ are strictly decreasing
$$
\hght(\b_1)>\hght(\b_2)>\cdots>\hght(\b_N)
$$

\begin{Def}
Let  $\l$ be a first layer quasi-dominant weight of co-length $N$ and 
$$
\l=\b_1+\cdots+\b_N
$$
its canonical expression. The height set and, respectively, the height vector of $\l$ are defined as 
$$
\Nght(\l):=\{\hght(\b_1), \dots, \hght(\b_N)\}
$$
and, respectively,
$$
\nght_\l:=(\hght(\b_1), \dots, \hght(\b_N))
$$
By definition $\Nght(0)$ and $\nght_0$ equal the empty set.
\end{Def}
For the above example $\Nght(\l)=\{3, 6, 7\}$. From the height set one can easily recover the height vector by writing the elements of the height set in decreasing order.

Since the largest possible height for a positive root is $n$ the height set can be regarded as a function 
$$
\Nght: Q^{(1)} \to \P([n])
$$

\begin{Prop}\label{heightsetbij}
 The map between the set of first layer quasi-dominant weights and set of subsets of $[n]$ which sends a weight to its height set is a bijection.
\end{Prop}
\begin{proof}
The map in the statement is of course well-defined since the largest possible height for a positive root is $n$. Let $S$ be an arbitrary subset of $[n]$. Our claim would follow if we show that there is a unique first layer quasi-dominant weight $\l$ such that  
\begin{equation}\label{eq9}
\Nght(\l)=S
\end{equation}

Denote by $N$ the cardinality of $S$. If $S$ is the empty set then it is  clear that the only possible weight satisfying \eqref{eq9} is the zero weight. Assume now that $N$ is positive and write $S=\{s_1,\dots,s_N\}$ with $s_1>\cdots>s_N$. 

First, remark that the sequence $$\{n+2-i-s_i\}_{1\leq i\leq N}$$ is weakly increasing. Second, 
\begin{equation}\label{eq10}
1\leq n+1-s_1\quad\text{and}\quad n+2-N-s_N\leq n+1-N
\end{equation}
The construction of $\l$ can be achieved as follows. We are looking for a first layer quasi-dominant weight of co-length $N$ such that if 
$$
\l=\b_1+\cdots+\b_N
$$
is its canonical expression then $\hght(\b_i)=s_i$ for all $1\leq i\leq N$. Since $\l$ is supposed to be quasi-dominant the coordinates which are equal to  $-1$ are on the rightmost possible positions. Hence, from the definition of the normal expression we know that for all $1\leq i\leq N$
$$
\b_i=\e_{\gamma(i)}-\e_{n+2-i}
$$
for some positive integer $\gamma(i)$. The condition on the height of $\b_i$ forces $$\gamma(i)=n+2-i-s_i$$
The inequalities \eqref{eq10} assure that $1\leq\gamma(i)\leq n+2-N$ and hence these positive integers do not interfere with the last $N$ coordinates. In consequence, we can define 
$$
\sum_{i=1}^N (e_{n+2-i-s_i}-e_{n+2-i})
$$
which is the desired  $\l$.
\end{proof}
Let us make explicit the definition of the inverse map
$$
\Nght^{-1}: \P([n]) \to Q^{(1),q+}
$$
The map sends the empty set to the zero weight. Otherwise, let 
$S=\{s_1,\dots,s_N\}$ be a subset of $[n]$ and assume that $s_1>\cdots>s_N$. Then, $$\Nght^{-1}(S)=(\l_1,\dots,\l_{n+1})$$  where   
\begin{subequations}
\begin{align}\label{eq11}
\l_k&:=|\{1\leq i\leq N~|~ n+2-i-s_i=k\}| \quad \text{for} \quad   1\leq k \leq n+1-N\\ \label{eq12}
\l_k&:=-1\quad \text{for}\quad  n+2-N\leq k\leq n+1 
\end{align}
\end{subequations}
The inclusion partial order on $\P([n])$ induces via this bijection a partial order on $Q^{(1),q+}$ which we also call inclusion and denote by $\subseteq$. More precisely, 
\begin{equation}
\l\subseteq \mu \quad \text{if and only if}\quad \Nght(\l)\subseteq \Nght(\mu)
\end{equation}

\begin{Def}\label{quasifundef} Let $\l$ be a first layer quasi-dominant weight. The function
$$
\QQ_\l:=\sum_{\mu\subseteq \l} \mm_\mu
$$
is called a (first layer) fundamental quasisymmetric function.
\end{Def}

It is clear from the definition that the fundamental quasisymmetric functions form a basis for $\QQ^{(1)}$. The characters of the first layer representations, being $S_{n+1}$ invariant are  elements of $\QQ^{(1)}$ and they can be  expressed as a linear combination of fundamental quasisymmetric functions. 

Let $\l$ be a first layer dominant weight. A quasi-dominant weight $\mu$ for which the coefficient of $\QQ_\mu$ in the expansion of the character $\chi_\l$ in the basis of fundamental quasisymmetric functions is non-zero is called a quasi-weight of $V_\l$. The set of quasi-weights of $V_\l$ will be denoted by ${\rm qwt}(\l)$. The coefficients in the expansion
\begin{equation}\label{fundamentalweights}
\chi_\l=\sum_{\mu\in{\rm qwt}(\l)} q_{\l\mu}\QQ_\mu
\end{equation}
will be called quasi-weight multiplicities.


\subsection{} 
We collect here a few remarks that will be used in the proof of the next Theorem. We introduce first some notation.

Assume that $n\geq 2$, $N\geq 1$ and let $\mu=(\mu_1,\dots,\mu_{n+1})$ be a first layer weight of co-length $N$  such that $\mu_1$ is  non-negative.  For any such element construct the following first layer weight  for the root system of type $A_{n-1}$ 
$$
 \mu^I:= (\mu_1+\mu_2,\mu_3,\dots,\mu_{n+1})  
$$
Comparing $\mu^I$ to $(0, \mu_1+\mu_2,\mu_3,\dots,\mu_{n+1})$ we observe that
\begin{equation}\label{height1}
\hght(\mu)=\hght(\mu^I)+\mu_1
\end{equation}
Remark also that $\mu^I$ has co-length $N$ unless $\mu_1>0$ and $\mu_2=-1$ in which case it has co-length $N-1$. Furthermore, if $\mu$ is quasi-dominant then $\mu^I$ is quasi-dominant.

Assume again that $n\geq 2$, $N\geq 1$ and this time let $\l=(\l_1,\dots,\l_{n+1})$ be a quasi-dominant first layer weight of co-length $N$ such that $\l_1$ is  positive.
Construct the following quasi-dominant first layer weight  of co-length $N-1$ for the root system of type $A_{n-1}$ 
$$
\l^{II}:= (\l_1-1,\l_2,\dots,\l_{n})
$$
Comparing $\l^{II}$ to $(\l_1-1,\l_2,\dots,\l_{n}, 0)$ we observe that
\begin{equation} \label{height2}
\hght(\l)=\hght(\l^{II})+n
\end{equation}

Let $$\l=\b_1+\cdots+\b_N $$ be the canonical decomposition of $\l$. As $\l_1>0$ it is clear that $\b_i=\e_1-\e_{n+2-i}$ for $1\leq i\leq \l_1$ and therefore
$$\nght(\l)=\nght(\l^I)+(1,\dots,1,0,\dots,0)\quad \text{and}\quad \nght(\l)=(n,\nght(\l^{II}))$$
In particular,
\begin{subequations}\label{factor1-2}
\begin{align}\label{factor1}
\frac{1-t^{-\nght_{\l^I}}}{1-t^{-\nght_\l}} &=\frac{1-t^{-n+\l_1}}{1-t^{-n}}\\ \label{factor2}
\frac{1-t^{-\nght_{\l^{II}}}}{1-t^{-\nght_\l}}&= \frac{1}{1-t^{-n}}
\end{align}
\end{subequations}

It is important to observe that if $w$ is a permutation such that $\ag_{w(\l)}<{\bf 0}$ and $w(\l)\in S_{n+1}\centerdot\l$ then $w(1)=1$. If $\l_2$ is non-negative then there are two possibilities: $w(2)=2$ or $w$ can be chosen such that  $w(n-N+2)=2$ (the first negative coordinate of $\l$ is on position $n-N+2$). Keeping this in mind and assuming that  $\l$ is quasi-dominant first layer weight of co-length $N$ such that $\l_1$ is  positive and $\l_2$ is non-negative, define
$$(S_{n+1}\centerdot\l)^{I}:=\{w(\l)\in S_{n+1}\centerdot\l~|~\ag_{w(\l)}<{\bf 0},~w(1)=1,~ w(2)=2\} $$ and $$(S_{n+1}\centerdot\l)^{II} :=\{w(\l)\in S_{n+1}\centerdot\l~|~\ag_{w(\l)}<{\bf 0},~w(1)=1,~ w(n+2-N)=2\}$$
The set 
 $$(S_{n+1}\centerdot\l)_{<{\bf 0}}:=\{\mu\in S_{n+1}\centerdot\l~|~\ag_\mu<{\bf 0} \}$$ is the disjoint union of 
$(S_{n+1}\centerdot\l)^{I} $ and $(S_{n+1}\centerdot\l)^{II}$.
\begin{Lm}\label{bij1} Let $\l$ be as above. 
The map $$(S_{n+1}\centerdot\l)^{I}\to (S_{n+1}\centerdot\l^I)_{<{\bf 0}}, \quad\mu \mapsto \mu^I$$
is a bijection. Moreover, for all elements of $(S_{n+1}\centerdot\l)^{I}$ we have
\begin{equation}\label{sum1}
c_\mu(t)=t^{\l_1}c_{\mu^I}(t)
\end{equation}
\end{Lm}
\begin{proof}  The first claim is clear from the definitions of the sets under consideration. The second claim follows from \eqref{height1} and from the fact that $\mu$ and ${\mu^I}$ have the same aggregate vectors.
\end{proof}
\begin{Lm}\label{bij2}  Let $\l$ be as above. 
The map $$(S_{n+1}\centerdot\l)^{II}\to (S_{n+1}\centerdot\l^{II})_{<{\bf 0}}, \quad\mu \mapsto \mu^I$$
is a bijection. Moreover, for all elements of $(S_{n+1}\centerdot\l)^{II}$ we have
\begin{equation}\label{sum2}
c_\mu(t)=t^{\l_1}(1-t^{-\l_1})c_{\mu^I}(t)
\end{equation}
\end{Lm}
\begin{proof} The first claim is again very easy to check. The second claim follows from \eqref{height1} and from the fact that
$$
(1-t^{\ag_\mu})=(1-t^{-\l_1})(1-t^{\ag_{\mu^{I}}})
$$
for all elements of $(S_{n+1}\centerdot\l)^{II}$.
\end{proof}

\subsection{} Our second  result is the following.

\begin{Thm}
 Let $\l$ be a first layer quasi-dominant weight. Then 
 \begin{equation}\label{monomial}
 \<1,\mm_\l\>_t=t^{\hght(\l)}(1-t^{-\nght_\l})
 \end{equation}
\end{Thm}
\begin{proof}  The statement is obviously satisfied by the zero weight. We can therefore assume  that $N\geq 1$. We will prove \eqref{monomial} by induction on $n\geq 1$.

If $n=1$ then $$\l=(1,-1)$$  and \eqref{monomial} is easily checked. 

Assume that $n\geq 2$.  First of all, if $\l_1=0$, observe that the map $$(S_{n+1}\centerdot\l)_{<{\bf 0}}\to (S_{n+1}\centerdot\l^I)_{<{\bf 0}}, \quad\mu \mapsto \mu^I$$
is a bijection and, 
$$
c_\mu(t)=c_{\mu^I}(t)
$$
for all elements of $(S_{n+1}\centerdot\l)_{<{\bf 0}}$. In consequence, from the induction hypothesis we obtain
\begin{eqnarray*}
\<1,\mm_\l\>_t&=&\<1,\mm_{\l^I}\>_t\\
&=&t^{\hght(\l^I)}(1-t^{-\nght_{\l^I}})\\
&=&t^{\hght(\l)}(1-t^{-\nght_\l})
\end{eqnarray*}
which is exactly our claim.

For the remainder of the argument, let us assume that $\l_1$ is positive. If $\l_2=-1$ then $n=N$ and $$\l=(N,-1,\dots,-1)$$
 It is clear that 
$$\ag_\l=-\nght_\l$$ and that $\l$ is the relevant element in $(S_{n+1}\centerdot\l)_{<{\bf 0}}$. Therefore,
\begin{eqnarray*}
\<1,\mm_\l\>_t&=&\<1,e^\l\>_t\\
&=&t^{\hght(\l)}(1-t^{-\nght_\l})
\end{eqnarray*}
which is exactly our claim. 

We can now focus our attention on the case when $\l_1$ is positive and  $\l_2$ is non-negative. The discussion leading to Lemma \ref{bij1} and Lemma \ref{bij2} applies, hence
\begin{eqnarray*}
\<1,\mm_\l\>_t&=&\sum_{\mu\in(S_{n+1}\centerdot\l)^{I} } c_\mu(t) +\sum_{\nu\in (S_{n+1}\centerdot\l)^{II}} c_\nu(t)\\
&=& t^{\l_1}\<1,\mm_{\l^I}\>_t+t^{\l_1}(1-t^{-\l_1})\<1,\mm_{\l^{II}}\>_t
\end{eqnarray*}
From the induction hypothesis $$\<1,\mm_{\l^I}\>_t= t^{\hght(\l^I)}(1-t^{-\nght_{\l^I}})\quad \text{and}\quad \<1,\mm_{\l^{II}}\>_t=t^{\hght(\l^{II})}(1-t^{-\nght_{\l^{II}}})$$
Keeping in mind the equations \eqref{height1},  \eqref{height2},  and \eqref{factor1-2} we obtain
\begin{eqnarray*}
\<1,\mm_\l\>_t&=& t^{\l_1}t^{\hght(\l^I)}(1-t^{-\nght_{\l^I}})+t^{\l_1}(1-t^{-\l_1})t^{\hght(\l^{II})}(1-t^{-\nght_{\l^{II}}})\\
&=& t^{\hght(\l)}(1-t^{-\nght_\l})\left( \frac{1-t^{-\nght_{\l^I}}}{1-t^{-\nght_\l}} +\frac{t^{-n+\l_1}(1-t^{-\l_1})(1-t^{-\nght_{\l^{II}}})}{1-t^{-\nght_\l}}\right)\\
&=& t^{\hght(\l)}(1-t^{-\nght_\l})\left(\frac{1-t^{-n+\l_1}}{1-t^{-n}} +\frac{t^{-n+\l_1}(1-t^{-\l_1})}{1-t^{-n}}\right)\\
&=& t^{\hght(\l)}(1-t^{-\nght_\l})
\end{eqnarray*}
which is exactly the desired formula. 
\end{proof}

\subsection{} We are now to prove our main result.  
\begin{Thm}\label{thmquasi}
 Let $\l$ be a first layer quasi-dominant weight. Then 
 \begin{equation}\label{quasi}
 \<1,\QQ_\l\>_t=t^{\hght(\l)}
 \end{equation}
\end{Thm}
\begin{proof} By the definition of fundamental quasisymmetric functions and the previous Theorem we obtain
 \begin{eqnarray*}
 \<1,\QQ_\l\>_t&=&\sum_{S\subseteq \Nght(\l)} (t^{S}-1)\\
 &=& \prod_{s\in \Nght(\l)}\left((t^{s}-1)+1\right)\\
 &=& t^{\hght(\l)}
 \end{eqnarray*}
 which is exactly our claim.
\end{proof}
We can immediately use Theorem \ref{quasi} to give a formula for generalized exponents.
\begin{Thm}\label{mainweakversion}
Let $\l$ be a first layer dominant weight. Then,
$$
E(V_\l)=\sum_{\mu\in{\rm qwt}(\l)} q_{\l\mu}t^{\hght(\mu)}
$$
\end{Thm}
\begin{proof} The result is a direct consequence of Theorem \ref{thmquasi} and \eqref{fundamentalweights}.
\end{proof}
At this point, an explicit description of  the set of quasi-weights and of the quasi-weight multiplicities is desirable. 
Some information  is known from the work of Gessel (in fact an immediate consequence of a very general result of Stanley on $P$-partitions) and it will be recalled in the next section. 


\section{Quasisymmetric polynomials}\label{quasipoly}

\subsection{}  The quasisymmetric functions introduced by Gessel \cite{gessel} are elements of the ring of power series in infinitely many variables. We are interested here only in the case when they have bounded degree and they depend on finitely many variables. To acknowledge this difference  and  to temporarily distinguish them from the quasisymmetric functions from the previous section they will be referred to  as {\sl quasisymmetric polynomials}.

Let  $f$ be an element of the polynomial ring $\Q[x_1,\dots, x_{n+1}]$. Monomials in this polynomial ring are written in the form $$x_{i_1}^{a_1}\cdots x_{i_k}^{a_k}$$ where $1\leq i_1<\cdots<i_k\leq n+1$ and the exponents $a_i$ are all positive integers. The polynomial $f$ is called quasisymmetric if the coefficient of any monomial $x_{i_1}^{a_1}\cdots x_{i_k}^{a_k}$ inside $f$ depends only on  $a_1,\dots, a_k$ (i.e. not on $1\leq i_1<\cdots<i_k\leq n+1$).
 
For any non-negative integer $k$ and any $k$-tuple ${\underline a}=(a_1,\dots, a_k)$ of positive integers define the quasisymmetric monomial 
$$\m_{\underline a}:=\sum_{1\leq i_1<\cdots<i_k\leq n+1} x_{i_1}^{a_1}\cdots x_{i_k}^{a_k}$$
The quasisymmetric monomials form a linear basis for the space of quasisymmetric polynomials.

If one is interested in studying homogeneous polynomials of some fixed degree the following notation is very useful. Fix $K$ a positive integer. The monomials $\m_{\underline a}$ for ${\underline a}=(a_1,\dots, a_k)$ such that $a_1+\cdots+a_k=K+1$ (the {\sl compositions} of $K+1$) form a linear basis of the space of quasisymmetric polynomials of degree $K+1$. The information about ${\underline a}$ can be encoded in the following subset of $[K]$
$$
A:=\{a_1,a_1+a_2,\dots, a_1+\cdots+a_{k-1}\}
$$
Conversely, the composition corresponding to the following subset of $[K]$
$$
B:=\{s_1,s_2,\dots, s_{l}\}
$$
where $s_1<\dots<s_{l}$ is 
\begin{equation}\label{eq15}
{\rm co}(S):=(s_1,s_2-s_1,\dots, s_l-s_{l-1}, K+1-s_l )
\end{equation}
Let us denote the set of compositions of $K+1$ by ${\rm Comp}(K+1)$. It is clear that the function
$$
{\rm co}: \P([K])\to {\rm Comp}(K+1)
$$
 is a bijection. We denote the inverse function by ${\rm co}^{-1}$.

The inclusion partial order relation on  $\P([K])$ induces via the above bijection a partial order relation on $ {\rm Comp}(K+1)$ that will be also referred to as inclusion and denoted by $\subseteq$. More precisely,
\begin{equation}
{\underline a} \subseteq {\underline b} \quad \text{if and only if}\quad {\rm co^{-1}}(\underline a)\subseteq {\rm co^{-1}}(\underline b)
\end{equation}

Another linear basis for the space of quasisymmetric polynomials of degree $K+1$ is the following
$$
\Qu_{\underline a}:=\sum_{\underline a\subseteq \underline b} \m_{\underline b}
$$
The polynomials $\Qu_{\underline a}$ are called {\sl fundamental quasisymmetric polynomials}. 


\subsection{}

A composition $\underline a = (a _{1}, \dots, a _{k})$  of $n+1$ is called partition if $a_1\geq\cdots\geq a_k$.
A partition could be represented its Young diagram which is consisting of  top-heavy, left-justified rows of square boxes, $a_i$ boxes in row $i$.
For example, the diagram of the partition $\underline a =(4,3,1)$ is
$$
 \tableau{{}&{}&{}&{}\\ {}&{}&{}\\{}}
 $$

For simplicity, we use the same symbol to denote a partition and its diagram.  A  semi-standard Young tableau  is a function
$\T: \underline a \rightarrow [n+1]$, which we picture as
assigning integer entries to the boxes of $\underline a $ such that the entries increase strictly  in columns and  increase weakly in rows.  Then,
\begin{equation}
x^{\T} := \prod _{u\in \underline a} x_{\T(u)},
\end{equation}
is a monomial of degree $n+1$ in $\Q[x_1,\dots, x_{n+1}]$. Denote the set of semi-standard Young tableaux of $\underline a$ by ${\rm SSYT}(\underline a)$.
If the tableau $\T$ is in addition a bijection  then $\T$ is called a standard Young tableau. The set of standard Young tableaux of $\underline a$ is denoted by ${\rm SYT}(\underline a)$.

Let $\underline a$ be a partition of $n+1$. One way to define the Schur function corresponding to $\underline a$ is as a generating function over the set of semi-standard Young tableaux of $\underline a$
$$
\Scal_{\underline a}:=\sum_{\T\in {\rm SSYT(\underline a)}} x^\T
$$

Let $\T$ be a standard Young tableau for the partition $\underline a$.  A descent of $\T$ is an integer $i$ such that $\T^{-1}(i+1)$ is a box in a lower row than $\T^{-1}(i)$. The set of descents of $\T$ is denoted by $\Des(\T)$. For example,  for 
$$
\T_1 = \tableau{1&2&6&8\\ 3&4&7\\ 5} \quad\quad\text{and}\quad\quad \T_2 = \tableau{1&2&4&6\\ 3&7&8\\ 5}
$$
we have $\Des(\T_1)=\Des(\T_2)=\{2, 4, 6\}$. The only result about quasisymmetric functions that we will need here is the following  particular case  of a result in \cite[pg. 295]{gessel}.
\begin{Thm} \label{schurexpansion}
Let $\underline a$ be a partition of $n+1$. Then,
\begin{equation}\label{schurexpansioneq}
\Scal_{\underline a}=\sum_{\T\in {\rm SYT(\underline a)}} \Qu_{\rm co(Des(\T))}
\end{equation}
\end{Thm}
Note that, as illustrated by the above example, the formula \eqref{schurexpansioneq} is not multiplicity free: different standard tableaux might have the same descent sets.

\subsection{} For the remainder of this section we explore the connection between homogeneous quasisymmetric polynomials of degree $n+1$ and first layer quasisymmetric functions.
First of all there is a straightforward connection between the indexing sets. Indeed, if $\l$ is a first layer quasi-dominant weight then $\l+\bf 1$ is a vector which has $\ell(\l)$ positive coordinates followed by  $\ell^*(\l)$ zero coordinates. Let us denote by $\underline{\l+\bf 1}$ the vector $\l+\bf 1$ truncated to the first $\ell(\l)$ coordinates. This is in fact a composition of $n+1$. Conversely, if $\underline a$ is a composition of $n+1$,  complete $\underline a$ to a $n+1$ vector by adding zero coordinates if necessary. We still use $\underline a$ to denote the outcome of this operation, using the words composition and vector to distinguish between the two objects. The vector $\underline a-\bf 1$ is a first layer quasi-dominant weight. The map 
$$
\varphi: Q^{(1),q+}\to {\rm Comp}(n+1),\quad \varphi(\l):=\underline{\l+\bf 1}
$$
is therefore a bijection. If we denote by $\Qu^{n+1}$ the space of quasisymmetric polynomials of degree $n+1$ then the linear map 
$$
\Phi: \QQ^{(1)}\to \Qu^{n+1}, \Phi(f):=(x_1\cdots x_{n+1})f
$$
is a well-defined isomorphism which sends $\mm_\l$ to $\m_{\varphi({\l})}$. The fundamental quasisymmetric functions and polynomials also correspond but for that we have to show that $\varphi$ is an anti-morphism of posets. Before we start presenting the argument we need to set up some notation. For a subset $A\subseteq [n]$, let
$$
(n+1)-A:=\{n+1-a~|~a\in A\}
$$
and denote by ${}^c\!A$ the complement of $A$ inside $[n]$. The function
$$
\phi:\P([n])\to\P([n]), \quad \phi(A):=(n+1)-{}^c\!A
$$
is an anti-involution of $\P([n])$ (i.e. reverses the inclusion partial order).

\subsection{} 

Let us start by recalling a couple of well-known fact about partitions. 
In what follows $1\leq N<n$ are two positive integers and $$\L=(\L_1,\dots, \L_{n-N}) $$ is a partition (eventually with zero parts) which fits inside the rectangle with $n-N$ rows and $N$ columns
$$
N\geq \L_1\geq \cdots \geq \L_{n-N} \geq 0
$$
To streamline some of the further considerations is it convenient to introduce
\begin{equation}\label{partitionconvention}
\L_0:=N\quad \text{and}\quad \L_{n-N+1}:=0
\end{equation}

Consider the partition dual to $\L$ to be the partition 
$$\L^\prime=(\L^\prime_1,\dots, \L^\prime_{N}) $$ defined by
$$
\L^\prime_i:=|\{1\leq j\leq n-N~|~ \L_l\geq i\}|
$$
for all $1\leq i\leq N$. The following  fact will be needed later
\begin{equation}\label{eq13}
\L_{i-1}-\L_i=|\{1\leq j\leq N~|~ \L^\prime_j=i-1\}|
\end{equation} 
for all $1\leq j\leq n-N+1$.
The second fact we need is a slight variation of a very basic fact which can be found for example in \cite[(1.7)]{macbook}.  The sets $$\L_h:=\{N+i-\L_i~|~ 1\leq i\leq n-N\}\quad \text{ and }\quad \L_v:=\{\L_j^\prime+N+1-j~|~ 1\leq j\leq N\}$$ are complementary subsets of $[n]$. The notation is motivated by the following interpretation of the above sets. Consider $\L$ being drawn inside the $(n-N)\times N$ rectangle and label the successive segments between $\L$ and its complement inside the $(n-N)\times N$ rectangle with the integers $1,\dots, n$, starting at the top right-hand corner and ending  at the lower left-hand corner of the rectangle.
  The elements of the set $\L_h$ are precisely the labels of the horizontal blocks and the elements of the set $\L_v$ are precisely the labels of the vertical blocks.

Let $0\leq N\leq n$ and let $$A=\{A_1,\dots, A_{n-N}\}$$ be a subset of $[n]$ of cardinality $n-N$, where $A_1<\cdots<A_{n-N}$. Denote also 
\begin{equation}\label{setconvention}
A_0:=0\quad \text{and}\quad A_{n-N+1}:=n+1
\end{equation}
Write $$
S:=\phi(A)=\{S_1,\dots, S_{N}\}
$$
where $S_1>\cdots>S_{N}$.

\begin{Lm} With the notation above
\begin{equation}\label{eq14}
A_k-A_{k-1}-1=|\{1\leq j\leq N~|~n+2-j-S_j=k\}|
\end{equation}
for all $1\leq k\leq n-N+1$.
\end{Lm}
\begin{proof}
Define $$
\L_i:=N-n-1+i+A_{n-N+1-i}
$$
for all $0\leq i\leq n-N+1$. It is straightforward to check that $\L:=(\L_1,\dots,\L_{n-N})$ is a partition  (eventually with zero parts) which fits inside the $(n-N)\times N$ rectangle
and that $\L_0=N$ and $\L_{n-N+1}=0$ in agreement with our convention \eqref{partitionconvention}.

Fix $1\leq k\leq n-N+1$. Then, keeping in mind \eqref{eq13} we obtain 
\begin{eqnarray*}
A_k-A_{k-1}-1&=&\L_{n-N+1-k}-\L_{n-N+1-(k-1)}\\
&=& |\{1\leq j\leq N~|~\L^\prime_j=n-N+1-k\}|
\end{eqnarray*}
However, $\L_h=(n+1)-A$ and since $\L_h$ and $\L_v$ are complementary inside $[n]$ we obtain that $\L_v=S$. In fact,
$$
S_j=\L^\prime_j+N+1-j
$$
for all $1\leq j\leq N$. Substituting this into the above formula for $A_k-A_{k-1}$ we obtain the desired statement.
\end{proof}
\begin{Lm} \label{lemma3}
With the notation above,
$$
\varphi(\Nght^{-1}(\phi(A)))={\rm co}(A)
$$
\end{Lm}
\begin{proof} Let $\l:=\Nght^{-1}(S)$ and $\underline a:={\rm co}(A)$. Since $S$ has size $N$, the first layer quasi-dominant weight $\l$ has co-length $N$ and therefore $\l_k=-1$ for $n-N+2\leq k\leq n+1$. Also $A$ has size $n-N$ and hence $\underline a$ is a composition of $n+1$ with $n-N+1$ parts. From \eqref{eq11}, \eqref{eq14}, and \eqref{eq15} we deduce that
$$
\l_k=a_k-1, \quad \text{for } 1\leq k\leq n-N+1 
$$
In consequence, $\underline{\l+\bf 1}=\underline a$, which is exactly our claim.
\end{proof}
\begin{Prop} The map 
$$
\varphi: (Q^{(1),q+},\subseteq)\to ({\rm Comp}(n+1), \subseteq)$$
is an anti-morphism of partially ordered sets. In particular the linear map 
$$
\Phi: \QQ^{(1)}\to \Qu^{n+1}
$$
sends $\QQ_\l$ to $\Qu_{\varphi(\l)}$ for all first layer quasi-dominant weights $\l$.
\end{Prop}
\begin{proof} Note that from Lemma \ref{lemma3} we know that
$$
\Nght^{-1}\circ\phi=\varphi^{-1}\circ{\rm co}
$$
Keeping in mind that $\phi$ is an involution we obtain
$$
\phi\circ\Nght={\rm co}^{-1}\circ\varphi
$$
Then,
\begin{eqnarray*}
\l\subseteq\mu & \text{iff}& \Nght(\l)\subseteq \Nght(\mu)\\
& \text{iff}& \phi(\Nght(\l))\supseteq \phi(\Nght(\mu))\\
& \text{iff}& {\rm co}^{-1}(\varphi(\l))\supseteq {\rm co}^{-1}(\varphi(\mu))\\
& \text{iff}& \varphi(\mu)\subseteq\varphi(\l)
\end{eqnarray*}
which shows that $\varphi$ is an anti-morphism of partially ordered sets. The remaining statement follows from the fact that $\Phi(\mm_\l)=\m_{\varphi(\l)}$ and the definition of fundamental quasisymmetric functions and, respectively, fundamental quasisymmetric polynomials.
\end{proof}
Combining this with Theorem \ref{schurexpansion} and the well-known fact that $$\Phi(\chi_\l)=\Scal_{\phi(\l)}$$ for all first layer dominant weights $\l$, we obtain the following result.

\begin{Prop}\label{quasi-characterexpansion}
Let $\l$ be a first layer dominant weight. Then,
$$
\chi_\l=\sum_{\T\in{\rm SYT(\varphi(\l))}} \QQ_{\phi(\Des(\T))}
$$
\end{Prop}
In particular, we have a description of the quasi-weights of $V_\l$
\begin{eqnarray*}
{\rm qwt}(\l)&=&\{\varphi^{-1}({\rm co}(\Des(\T)))~|~\T\in{\rm SYT}(\l)\}\\
&=&\{\Nght^{-1}({\phi}(\Des(\T)))~|~\T\in{\rm SYT}(\l)\}
\end{eqnarray*}
It would be interesting to describe the set of quasi-weights abstractly, without reference to above construction. For example, it is clear that the maximal elements of $\wt(\l)\cap Q^{(1),q+}$ with respect to the inclusion partial order must be in ${\rm qwt}(\l)$.

If $\mu=\Nght^{-1}({\phi}(\Des(\T)))$ for some $\T\in{\rm SYT}(\l)$, then 
\begin{eqnarray*}
\hght(\mu)&=&\sum_{b\in\Nght(\mu)} b\\
&=&\sum_{b\in{\phi}(\Des(\T))} b\\
&=&\sum_{a\in{}^c\!\Des(\T)} (n+1-a)
\end{eqnarray*}
Motivated by this fact, for any standard Young tableau $\T$, we define the positive integer
\begin{equation}\label{tableauheight}
\hght(\T):=\sum_{a\in{}^c\!\Des(\T)} (n+1-a)
\end{equation}
which we call the height of $\T$.
\subsection{}\label{maintheorem} We are now ready to prove our last result. 

\begin{Thm}\label{main}
Let $\l$ be a first layer dominant weight. Then,
\begin{equation}\label{mainformula}
E(V_\l)=\sum_{\T\in{\rm SYT}(\varphi(\l))} t^{\hght(\T)}
\end{equation}
\end{Thm}
\begin{proof} Straightforward from Theorem \ref{thmquasi}, Proposition \ref{quasi-characterexpansion} and \eqref{tableauheight}.
\end{proof}

\section{Comparison with the charge formula} \label{lascoux}
	
 The polynomials $E(V_\l)$ are in fact particular cases of the celebrated Kostka-Foulkes polynomials $K_{\l\mu}(t)$ which, in one formulation,  are the entries of the change of basis matrix from Schur functions to Hall-Littlewood symmetric functions. Indeed,
\begin{equation}\label{kostka}
E(V_\l)=K_{\varphi(\l),\varphi(0)}(t)
\end{equation}
Lascoux and Sch\" utzenberger \cite{lascoux} discovered a truly remarkable non-negative formula for the Kostka-Foulkes polynomials, showing that they are in fact the generating function of a tableaux statistic that they called charge.  We briefly recall the definition here, not in  full generality but in the context that is relevant for us. We refer the reader to 
\cite{lascoux}, \cite[Section 2.4]{butler}, and \cite[III, 6, pg. 242]{macbook} for the most general statements.

Let $w$ be an element of $S_{n+1}$. Identify $w$ with the sequence $w(1),w(2),\dots,w(n+1)$. Each integer $i\geq 2$ appearing in  this sequence has a contribution ${\rm chc}(i)$  to the charge of $w$ as follows: ${\rm chc}(i)$ is either zero or $n+1-(i-1)$ depending on whether $i$  appears in $w$ to the right or to the left of $i-1$. The charge of $w$ is $${\rm ch}(w):={\rm chc}(2)+\cdots+{\rm chc}(n+1) $$

If $\T$ is a standard Young tableau, read the entries of $\T$ from right to left in consecutive rows starting from the top, to obtain a permutation denoted $w(\T)$. The charge of $\T$ is defined to be the charge of $w(\T)$. The Lascoux--Sch\" utzenberger formula is
\begin{equation}\label{charge}
K_{\varphi(\l),\varphi(0)}(t)=\sum_{\T\in{\rm SYT}(\varphi(\l))} t^{{\rm ch}(\T)}
\end{equation}
However, this is precisely \eqref{mainformula} since \begin{equation}\label{chargeisheight}\hght(\T)={\rm ch}(\T)\end{equation} for all standard Young tableaux. 

Indeed, assume that $\T$ is a standard Young tableaux and let $2\leq i\leq n+1$. If $i-1$ is a descent of $\T$, then $i$ appears in $w(\T)$ to the right of $i-1$ and therefore ${\rm chc}(i)$ is zero. If $i-1$ is not a descent of $\T$, then $i$ appears in $\T$ immediately to the right of  $i-1$ and hence in $w(\T)$ to the left of $i-1$. Therefore, if $i-1$ is not a descent of $\T$, then ${\rm chc}(i)$ equals $n+1-(i-1)$. By the definition of charge,
$$
{\rm ch}(\T)=\sum_{j\in{}^c\!\Des(\T)} (n+1-j)
$$
which is exactly the definition  \eqref{tableauheight} of $\hght(\T)$.


\end{document}